\documentclass[a4paper]{amsart} 
\usepackage[a4paper]{geometry}
\usepackage[utf8]{inputenc}
\usepackage{babel}
\usepackage{amsthm,amsmath,amssymb,amsfonts,amscd}
\usepackage{mathtools}
\usepackage{hyperref}
\usepackage{xcolor}
\usepackage{subfiles}
\usepackage{enumerate}
\usepackage{comment}
\usepackage{graphicx} 
\usepackage{thm-restate}
\usepackage{esint}

\newtheorem{theorem}{Theorem}[section]
\newtheorem{lemma}[theorem]{Lemma}
\newtheorem{prop}[theorem]{Proposition}
\theoremstyle{definition}

\theoremstyle{remark}
\newtheorem{remark}{Remark}

\newcommand{\R}{\ensuremath{\mathbb{R}}}
\newcommand{\T}{\ensuremath{\mathbb{T}}}
\newcommand{\vag}{v_{\mathrm{Ag}}}
\newcommand{\kag}{K_{\mathrm{Ag}}}
\newcommand{\vbs}{v_{\mathrm{BS}}}

\newcommand{\pnorm}[2]{\|#2\|_{#1}}
\newcommand{\diff}{\mathrm{d}}

\usepackage{soul}

\usepackage[disable,textsize=tiny,textwidth=27mm]{todonotes}

\author{Juan Carlos Cantero}
\address{Departament de Matemàtiques -- Universitat Autònoma de Barcelona\\
08193 Bellaterra, Barcelona, Spain}

\email{juancarlos.cantero@uab.cat}

\author{Marc Magaña}
\address{Departament de Matemàtiques -- Universitat Autònoma de Barcelona\\
08193 Bellaterra, Barcelona, Spain}
\email{marc.magana@uab.cat}

\title[$C^2$ regularity of planar patches]{%
  Persistence of $C^2$ regularity for aggregation patches
  and a classification for linear transforms of the Newtonian gradient}

\date{\today}

\begin{document}

\maketitle

\begin{abstract}
    We prove persistence of $C^2$ boundary regularity for patch solutions of the two-dimensional aggregation transport equation, in contrast with the known $C^2$ ill-posedness of Euler vortex patches. More generally, we consider transport equations with kernels of the form $L(\nabla N)$, where $N$ is the fundamental solution of the Laplacian and $L$ is a linear map. We show that $C^2$ regularity persists for every bounded simply connected initial patch whose boundary is a regular $C^2$ Jordan curve if and only if $L$ is a scalar multiple of the identity. Our approach builds on the curvature formulation introduced by Kiselev and Luo.
\end{abstract}

\tableofcontents

\newpage

\section{Introduction} 
Active scalar equations form an important class of nonlocal transport equations arising in fluid dynamics and mathematical biology. In these equations, a scalar quantity is transported by a velocity field that depends nonlocally on the scalar itself. Patch solutions provide a natural setting in which to study the propagation of geometric regularity for this class of equations. In this paper we study this question for the two-dimensional aggregation transport equation and, more generally, for kernels obtained by applying a linear map to the gradient of the Newtonian potential. Our principal concern is the persistence of $C^2$ regularity for patch boundaries. 

We consider the Cauchy problem 
\begin{equation}\label{transportequation}
\begin{cases}    \partial_t\rho+v\cdot\nabla\rho=0, \\
 v(\cdot,t)=K*\rho(\cdot,t), \\
 \rho|_{t=0}=\chi_{D_0},
\end{cases}
\end{equation}
where $D_0\subset\mathbb R^2$ is a bounded domain and $K$ is the interaction kernel. For sufficiently regular velocity fields, the transport equation preserves the value of the scalar along the trajectories of the flow. Hence, if $\rho(\cdot,0)=\chi_{D_0}$, then $\rho(\cdot,t)=\chi_{D_t}$, where $D_t$ is a bounded domain evolving under the flow generated by $v$. Such solutions are called patch solutions.

Two important examples of this equation are the two-dimensional Euler equation and the aggregation equation. For Euler, the scalar $\rho$ is the vorticity and it is usually denoted by $\omega$, and the velocity is given by the Biot–Savart law, whose kernel is
$$
K_{\mathrm{BS}}(x)=\nabla^\perp N(x)=\frac{x^\perp}{2\pi|x|^2},
\qquad x^\perp=(-x_2,x_1),
$$
where $N(x)=\frac{1}{2\pi}\log|x|$ is the fundamental solution of the Laplacian in two dimensions. Thus, the Euler velocity satisfies
$\vbs=\nabla^\perp N*\rho
=-\nabla^\perp(-\Delta)^{-1}\rho,$ $
\operatorname{div}v_{\mathrm{BS}}=0.
$

The aggregation equation considered here corresponds to the kernel
\begin{equation}\label{eq:aggregation-kernel}
K_{\mathrm{Ag}}(x)=-\nabla N(x)
=-\frac{x}{2\pi|x|^2}=K_{\mathrm{BS}}(x)^\perp.
\end{equation}
Consequently, the aggregation velocity satisfies $v_{\mathrm{Ag}}=\nabla(-\Delta)^{-1}\rho$, and $\operatorname{div}v_{\mathrm{Ag}}=-\rho$ in the sense of distributions. 

The literature on vortex patches for the two-dimensional Euler equation is extensive, and we highlight here the results most relevant to our study.

At the level of weak solutions, Yudovich's theorem (see \cite{Y}) ensures global existence and uniqueness for initial vorticity data in $L^1\cap L^{\infty}.$ In particular, this provides global well-posedness for vortex patches. However, this result does not provide control of the regularity of the boundary of the patch. This question was resolved in the Hölder setting by Chemin \cite{Ch}, who proved global persistence of $C^{1,\alpha}$ regularity (\(0<\alpha<1\)) using paradifferential calculus. Shortly thereafter, Bertozzi and Constantin \cite{BC} obtained an alternative proof based on geometric methods in classical analysis. Their arguments can be extended to show persistence of $C^{k,\alpha}$ regularity, for any integer $k\ge 1$ and  $0<\alpha<1$. The endpoint regularity question for $C^2$ and $C^{1,1}$ was settled by Kiselev and Luo in \cite{KL}. They constructed \(C^2\) initial patches whose boundary curvature is unbounded for every sufficiently small positive time, thereby establishing ill-posedness in both \(C^2\) and \(C^{1,1}\).

For the specific case of rotating vortex patches, several questions concerning their existence and regularity have been addressed in the literature. The Rankine vortex \cite{R1858}, corresponding to a circular patch, and Kirchhoff's rotating ellipse \cite{K1876} are classical examples of such solutions, with explicitly known shapes. Deem and Zabusky \cite{DZ1978} numerically computed the first families of non-trivial rotating patches, the so-called V-states, emanating from the disk, and Burbea \cite{B} later gave a rigorous construction of these branches by bifurcation from the disk. Hmidi, Mateu and Verdera \cite{HMV} proved smoothness of these patches, while Hmidi and Mateu \cite{HM} studied bifurcation from Kirchhoff ellipses. Later on, spatial analyticity was established for such solutions. Castro, C\'ordoba and G\'omez-Serrano \cite{CCGS} constructed uniformly rotating analytic patches. More recently, Castro-L\'opez and G\'omez-Serrano \cite{CLGS} established the existence of uniformly rotating analytic patches with non-convex boundaries. 

Despite its importance in applications such as biological swarming and the evolution of vortex densities in superconductors (see \cite{BLL} and the references therein), the mathematical literature on the aggregation equation is considerably more limited. The Newtonian aggregation equation is usually written in divergence form,
\begin{equation}\label{aggregationdivergence}
\begin{cases}
\partial_s\rho+\operatorname{div}(\rho v)=0,\\
v(\cdot,s)=K_{\mathrm{Ag}}*\rho(\cdot,s),\\
\rho(\cdot,0)=\rho_0.
\end{cases}
\end{equation}
In a framework analogous to Yudovich's theory for the Euler equation, Bertozzi, Laurent and L\'eger \cite{BLL} proved that, for compactly supported $\rho_0\in L^1\cap L^\infty$, \eqref{aggregationdivergence} has a unique distributional solution in $C([0,T];L^1)\cap L^\infty$ for every $T<(\sup\rho_0)^{-1}$. For $\rho_0=\chi_{D_0}$ this solution is
\begin{equation}\label{eq:timechange}
\rho(x,s)=\frac{1}{1-s}\,\chi_{D_t}(x),\qquad t=-\log(1-s),\quad 0\le s<1,
\end{equation} where $\chi_{D_t}$ is the solution of \eqref{transportequation} with $K=K_{\mathrm{Ag}}$ considered in Theorem~\ref{thm:C2persistence} (see \cite[Section~4]{BLL}). In the divergence formulation the mass is conserved while, by \eqref{eq:timechange}, the area of the patch $|D_{t}|=(1-s)|D_0|$ tends to zero as $s\to1$, so the density blows up at the finite time $s=1$. The patch is expected to concentrate on a lower-dimensional set, called the skeleton: a disk collapses to its center and an ellipse onto a segment, while for more general patches the numerical simulations in \cite{BLL} show branched structures. Rigorous descriptions of the limiting set are available only in particular cases, such as one-fold symmetric patches \cite{HL} and the one-dimensional problem \cite{CO}. In the transport formulation, instead, the patch is defined for all $t\ge0$, since the collapse time corresponds to $t=\infty$. Bertozzi, Garnett, Laurent and Verdera \cite{BGLV} proved, in every dimension, that $C^{1,\alpha}$ regularity of the boundary of aggregation patches persists for all times.

Thus, both the two-dimensional Euler equation and the aggregation equation preserve $C^{1,\alpha}$ regularity of patch boundaries. A natural question was whether they should have the same behavior for $C^2$ boundaries. The answer is negative. We prove that the two-dimensional aggregation transport equation preserves $C^2$ boundaries. This contrast motivates our main result.

\begin{restatable}[Persistence of $C^2$ patches]{theorem}{persistence}
\label{thm:C2persistence}
Let $D_0\subset\mathbb R^2$ be a bounded simply connected domain
whose boundary is a regular $C^2$ Jordan curve. Let
$\rho(x,t)=\chi_{D_t}(x)$ be the unique solution of the aggregation transport equation
\begin{equation*}
\begin{cases} \partial_t\rho+\vag\cdot\nabla\rho=0, \\
 \vag=-\nabla N*\rho,
 \qquad  \qquad
 N(x)=\frac1{2\pi}\log|x|,\\
 \rho(\cdot,0)=\chi_{D_0}.
\end{cases}
\end{equation*}
Then $\partial D_t$ is a regular $C^2$ Jordan curve for every
$t\ge0$.
\end{restatable}

Both Euler and aggregation kernels belong to a broader family of kernels of the form $K_L=L(\nabla N),$ where $L:\mathbb R^2\to\mathbb R^2$ is a linear map. Another notable member of this family corresponds to \(L=\begin{pmatrix}1&0\\ 0&-1\end{pmatrix}\), and is given in complex notation by the Cauchy kernel \(K(z)=\frac{1}{2\pi z}\). Although the associated velocity field is not divergence free, the existence of Yudovich-type solutions has been established for this kernel by Clop and Sengupta \cite{CS} via quasiconformal map techniques, and their uniqueness has recently been proved by Cantero, Clop, Faraco, Horas, Mateu, Orobitg and Sengupta \cite{CCFHMOS}.

Cantero, Mateu, Orobitg and Verdera \cite{CMOV} established a more general patch persistence result in arbitrary dimension $d$, replacing the Biot-Savart and aggregation kernels with a general $K$ that is odd, homogeneous of degree $-(d-1)$ and $C^2$ away from the origin. Under these assumptions, they proved global persistence of $C^{1,\alpha}$ patch boundaries, together with uniqueness in the corresponding class of regular patches. Their theorem includes Euler and aggregation, but also kernels which are neither divergence free nor gradients of scalar potentials. Later, Verdera~\cite{V} gave a different planar proof for their family of kernels. In particular, both results apply in the plane to kernels of the form $K_L(x)=L(\nabla N(x))$. 

This raises the question of whether \(C^2\) regularity is also preserved for the family of kernels \(K_L\). Our next theorem answers this question by characterizing precisely the linear maps \(L\) for which persistence of \(C^2\) regularity holds. In particular, it shows that the aggregation kernel is exceptional within this family: persistence for every \(C^2\) initial patch holds precisely when \(L\) is a scalar multiple of the identity map.

\begin{restatable}{theorem}{classification}
\label{thm:C2-classification}
Let $L:\mathbb R^2\to\mathbb R^2$ be linear, and consider equation \eqref{transportequation} with kernel $K=K_L,\; K_L(x)=L(\nabla N(x)).$ Then, $C^2$ regularity  of patch boundaries is preserved for every bounded simply connected initial domain whose boundary is a regular $C^2$ Jordan curve if and only if $L=c\,\mathrm{Id}$ for some $c\in\mathbb R$.
More precisely, if $L$ is not of this form, there exist such a domain $D_0$ and $T>0$ such that the curvature of its
corresponding patch solution satisfies $  \kappa(t)\notin L^\infty
    \text{ for every }0<t<T.$
\end{restatable}

Both theorems also hold with $C^{1,1}$ in place of $C^2$. In fact, for $L=c\,\mathrm{Id}$ singularities of the initial curvature are only transported and multiplied by a continuous positive factor (Remarks~\ref{rmk:C11} and~\ref{rmk:C11L}).

The persistence of \(C^2\) regularity for aggregation patches is particularly striking in view of the contrasting behavior of Euler patches established by Kiselev and Luo. Indeed, within the family \(K_L\), the Euler kernel, up to scalar multiples, is the only nontrivial one that generates a divergence-free velocity field, yet \(C^2\) regularity may fail to persist for Euler patches. By contrast, the aggregation velocity satisfies \(\operatorname{div}v_{\mathrm{Ag}}=-\rho\), which is discontinuous across the patch boundary. Thus, the persistence of $C^2$ regularity for aggregation cannot be attributed to incompressibility. The reason lies in the evolution equation for the curvature $\kappa$ of the boundary. By \eqref{evolcurv}, the shape of the boundary changes through the normal component of $\partial_s^2v$, the second arc-length derivative of the velocity. In $\partial_s^2v$, the curvature at nearby points of the boundary enters through a singular integral. Since these points lie almost on the tangent line, they are seen in the directions $\pm\tau$, where the kernel takes the values $K_L(\pm\tau)=\pm\frac{1}{2\pi}L\tau$; hence, up to smoother terms, this singular integral is a multiple of the vector $L\tau$ times the Hilbert transform $\mathcal H\kappa$. Its normal component produces in the curvature equation the term $q_L\,\mathcal H\kappa$, with
\[
q_L=\tfrac12(L\tau)\cdot n=\pi\,K_L(\tau)\cdot n
\]
(Proposition~\ref{prop:Lsystem}). Since the Hilbert transform of a continuous function need not be bounded, this term destroys $C^2$ regularity for suitable initial data, as in \cite{KL} for Euler, where $q_L\equiv-\frac12$. For $L=c\,\mathrm{Id}$, instead, $L\tau=c\tau$ is tangent to the boundary, so the singular term only moves points along the curve and does not change its shape. This happens for every tangent direction exactly when $K_L(x)$ is parallel to $x$ for every $x$, that is, when $L=c\,\mathrm{Id}$ (equivalently, when the kernel $K_L=\nabla(cN)$ is a gradient).

\subsection{Outline of the paper}
In Section~\ref{sec:w2p} we prove the $W^{2,p}$ regularity. In Section~\ref{sec:curvature} we derive the curvature equation for aggregation. Section~\ref{sec:C2agg} contains the main theorem: persistence of $C^2$ regularity for the aggregation model. Finally, Section~\ref{sec:generalK} deals with general kernels of the form $L(\nabla N)$.

\subsection{Notation, conventions, useful relations and function spaces} 

To make the paper easier to read, we begin by introducing the notation and conventions that will be used throughout the paper. For the reader's convenience, we also recall some useful geometric identities and define the function spaces in which our arguments will be carried out.

For a vector \(v=(v_1,v_2)\in\mathbb R^2\), we define its perpendicular by \(v^\perp=(-v_2,v_1)\). We identify \(\mathbb T\) with \(\mathbb R/(2\pi\mathbb Z)\). For a curve \(\T\ni\xi\mapsto\gamma(\xi)\), we define the arc-length metric \(g(\xi)=|\partial_\xi\gamma(\xi)|\), the unit tangent vector $\tau(\xi)=\frac{\partial_\xi\gamma(\xi)}{g(\xi)},$ and the normal vector \(n(\xi)=-\tau^\perp(\xi)\). If a closed simple curve is parametrized counterclockwise, then \(n(\xi)\) is the outer normal vector with this convention.

The arc-length parameter is defined by $s=\int_0^\xi g(\eta)\,\mathrm \diff \eta.$ The chain rule gives the relationship \(\partial_\xi=g\partial_s\). Moreover, the tangent and normal vectors satisfy the Frenet formulas
$$
\partial_s\tau=-\kappa n,\qquad
\partial_s n=\kappa\tau,
$$
where $\kappa$ is the signed curvature of $\gamma$.

Regarding the function spaces we will work with, we set the quantities 
$$\Gamma(\gamma):= \sup_{\xi,\eta \in \mathbb{T},\xi \neq \eta} \frac{\left|\xi-\eta\right|}{\left|\gamma(\xi)-\gamma(\eta)\right|},\quad \left|g\right|_*=\max\{\pnorm{\infty}{g},\pnorm{\infty}{1/g}\},$$

where \(|\xi-\eta|\) denotes the periodic distance on \(\mathbb T\). The quantity \(\Gamma(\gamma)\) is the inverse arc-chord ratio and controls self-intersections and degeneration of the parametrization, while \(|g|_*\) controls the size of the parametrization speed and its reciprocal. Note that \(|g|_*\) is not a norm.

Then we will work in
$$X_p:=\{\gamma: \mathbb{T}\to \mathbb{R}^2: \gamma \in W^{2,p}(\mathbb{T}) \text{ with } \Gamma(\gamma) < \infty\},$$
which is an open subset of the Banach space $W^{2,p}(\mathbb T;\mathbb R^2)$, and in a closed subset of it,
$$B_p^M:=\{\gamma \in X_p: \|\gamma\|_{W^{2,p}(\mathbb{T})}\le M,\ |g|_*\le M,\ \Gamma(\gamma) \le M\}.$$
When a parametrization represents the boundary of a positive patch, we take it counterclockwise.

\section{\texorpdfstring{$W^{2,p}$ regularity}{W2p regularity}}\label{sec:w2p}
To study the persistence of \(W^{2,p}\) regularity of the boundary of an aggregation patch, we use a contour dynamics formulation. This formulation describes the evolution of the boundary through an equation for its parametrization, reducing the patch evolution problem to a nonlocal evolution equation for a curve. The contour dynamics formulation is classical for the two-dimensional Euler equations. It was introduced by Deem and Zabusky \cite{DZ1978} and Zabusky, Hughes and Roberts \cite{ZHR}; see also Majda \cite{Majda1986} and Majda and Bertozzi \cite{MajdaBertozzi} for background and a detailed treatment. For the aggregation equation, and for the kernels of Section \ref{sec:generalK}, the contour dynamics equation is obtained analogously; see \cite[Section 2]{BGLV} and \cite[Section 2.1]{CMOV}.

Let \(\rho(\cdot,t)=\chi_{D_t}\) be the aggregation patch, and let \(\gamma(\cdot,t)\) be a counterclockwise parametrization of $\partial D_t$ whose points move with the velocity field, $\partial_t\gamma(\xi,t)=v_{\mathrm{Ag}}(\gamma(\xi,t),t)$. Since \(v_{\mathrm{Ag}}=-\nabla N*\rho\), where \(N(x)=\frac{1}{2\pi}\log|x|\), Green's formula yields the contour dynamics equation (CDE)
\begin{equation}
\partial_t\gamma(\xi,t)=
-\frac{1}{2\pi}
\int_{\mathbb T}
\log\bigl|\gamma(\xi,t)-\gamma(\eta,t)\bigr|
(\partial_{\eta}\gamma(\eta,t))^\perp\, \diff \eta.
\tag{CDE}\label{CDE}
\end{equation}

Local existence and uniqueness in $W^{2,p}$ follow from a standard fixed-point argument, which yields a unique solution on a finite time interval. We omit the technical details and refer the interested reader to \cite[Section 2]{KL}, where the same argument is developed for vortex patches. This is possible due to the fact that the velocity field considered here is obtained by rotating the Euler velocity field, as we argued before. Nevertheless, we include the formulas for the velocity derivatives, which are needed throughout. More precisely, the  Euler and aggregation velocities at $\gamma(\xi,t)$ are
\begin{align*}
v_{\mathrm{Ag}}(\xi)
&=\frac1{2\pi}\int_{\mathbb T}n(\eta)\log|\gamma(\xi)-\gamma(\eta)|\,g(\eta)\,\diff\eta,\\
v_{\mathrm{BS}}(\xi)
&=-\frac1{2\pi}\int_{\mathbb T}\tau(\eta)\log|\gamma(\xi)-\gamma(\eta)|\,g(\eta)\,\diff\eta.
\end{align*}

We observe that the only difference is the factor $n(\eta)$ instead of the factor $-\tau(\eta)$ in the integrand. Indeed, both kernels are of the form $L(\nabla N),$ with $L=\begin{pmatrix}  0 & -1 \\ 1 & 0 \end{pmatrix}$ for Euler and $L=-\operatorname{Id}$ for aggregation, and the velocity at $\gamma(\xi)$ is then $-\frac{1}{2\pi}\int_{\mathbb T} Ln(\eta)\log|\gamma(\xi)-\gamma(\eta)|\,g(\eta)\,\diff \eta.$ Since $L$ is constant, it commutes with $\partial_\xi$ and $\partial_s,$ so all the computations involving derivatives in $\xi$ hold for every constant matrix $L$. Since the aggregation velocity is obtained by applying an orthogonal rotation to the Euler velocity, the corresponding formulas are obtained by applying the same rotation to the vector-valued factors in the formulas for Euler. In particular, whenever they appear outside the scalar kernels, $\tau(\eta)$ is replaced by $-n(\eta)$, while $n(\eta)$ is replaced by $\tau(\eta)$. Therefore,
\cite[Proposition 2.4]{KL} yields the following result for aggregation
patches. 
\begin{prop}[Arc-length derivatives of the velocity at $\partial D$]
    Let $1<p<\infty$, $M>1$, and let $\gamma\in B_p^M$ be a counterclockwise parametrization of the boundary of a domain $D$. Then $v_{\mathrm{Ag}}\in W^{2,p}(\mathbb T;\mathbb R^2)$. Its first and second arc-length derivatives are given almost everywhere by the following formulas.

\begin{align}
\label{firstderofv}\partial_s v_{\mathrm{Ag}}(\xi)=&
\frac1{2\pi}\operatorname{p.v.}
\int_{\mathbb T}
n(\eta)
\frac{(\gamma(\xi)-\gamma(\eta))\cdot\tau(\xi)}
{|\gamma(\xi)-\gamma(\eta)|^2}
\,g(\eta)\,\diff\eta,\\
\label{secondderofv}\partial_s^2v_{\mathrm{Ag}}(\xi)
={}&
\frac1{2\pi}\operatorname{p.v.}
\int_{\mathbb T}
\kappa(\eta)\tau(\eta)
\frac{(\gamma(\xi)-\gamma(\eta))\cdot\tau(\xi)}
{|\gamma(\xi)-\gamma(\eta)|^2}
\,g(\eta)\,\diff\eta
\\[1mm]
&+\frac1{2\pi}\int_{\mathbb T}
n(\eta)
\frac{
(\tau(\xi)-\tau(\eta))\cdot\tau(\xi)
-\kappa(\xi)(\gamma(\xi)-\gamma(\eta))\cdot n(\xi)
}{
|\gamma(\xi)-\gamma(\eta)|^2
}
\,g(\eta)\,\diff\eta  \notag
\\[1mm]
&-\frac1{\pi}\int_{\mathbb T}
n(\eta)
\frac{
\bigl((\gamma(\xi)-\gamma(\eta))\cdot\tau(\xi)\bigr)
\bigl((\gamma(\xi)-\gamma(\eta))\cdot(\tau(\xi)-\tau(\eta))\bigr)
}{
|\gamma(\xi)-\gamma(\eta)|^4
}
\,g(\eta)\,\diff\eta,\notag
\end{align}
where the right-hand sides of \eqref{firstderofv} and \eqref{secondderofv} are well-defined in \(W^{1,p}(\mathbb T)\) and \(L^p(\mathbb T)\), respectively, and agree with the corresponding weak arc-length derivatives. Furthermore, 
$$
\|v_{\mathrm{Ag}}\|_{W^{2,p}(\mathbb T)}
\le C_{p,M}.$$
\end{prop}

Now, since the kernels involved in \eqref{firstderofv} and \eqref{secondderofv} are of the same type as those for the Euler case, we can recover (see \cite[Propositions 2.7 and 2.8]{KL}) the following bounds for $\gamma_1,\gamma_2\in B_p^M,$
$$\pnorm{L^p(\mathbb{T})}{\partial^2_s v(\gamma_1)-\partial_s^2 v(\gamma_2)},\pnorm{L^{\infty}(\mathbb{T})}{\vag(\gamma_1)-\vag(\gamma_2)} \le C_{M}\pnorm{W^{2,p}}{\gamma_1-\gamma_2}.$$

Let $\gamma_0$ be a parametrization of $\partial D_0$ with finite inverse arc-chord ratio. The fixed-point argument in \cite[Proposition 2.9]{KL} gives the local-in-time existence and uniqueness in $W^{2,p}$. Consequently, we have the following theorem.

\begin{theorem}[Local-in-time existence and uniqueness in $W^{2,p}$]
\label{thm:localw2p}
    Let $1<p<\infty$. If $\gamma_0 \in X_p$ then there is a unique local-in-time solution $\gamma$ to equation \eqref{CDE} in $C([0,T];X_p)$ for some $T>0$ with initial condition $\gamma(0)=\gamma_0.$
\end{theorem}

The scheme to prove global-in-time regularity of solutions is not that classical. We follow the clever argument in \cite[Section 3]{KL}. The idea is to embed $W^{2,p}$ in a space where global solutions exist and then get an improved $W^{2,p}$ estimate taking advantage of the fact that the solution is $C^{1,\alpha}$. More specifically, since $W^{2,p}(\mathbb{T}) \hookrightarrow C^{1,\alpha}(\mathbb{T})$ provided $\alpha=1-1/p$, the global $C^{1,\alpha}$ theory controls all the quantities defining $B^M_p$ except $\|\gamma\|_{W^{2,p}}$, which is then controlled by an improved estimate. The only ingredient that depends on the equation is \cite[Proposition~3.1]{KL}: the solution of \eqref{CDE} continues for all times in $C^{1,\alpha}$ and, for every $T>0$,
\begin{equation}\label{eq:C1alpha}
\sup_{0\le t\le T}\Bigl(|g(t)|_*+\Gamma(\gamma(t))+\|\tau(t)\|_{C^{\alpha}(\mathbb T)}+\|\partial_s v(\gamma(t))\|_{L^\infty(\mathbb T)}\Bigr)<\infty .
\end{equation}

For aggregation patches, by the uniqueness in \cite[Theorem~2.1]{BGLV} the local solution in Theorem~\ref{thm:localw2p} coincides with the global $C^{1,\alpha}$ solution of \cite[Sections~2 and~7]{BGLV}, and \eqref{eq:C1alpha} follows from the a priori bounds in \cite[(7.15)--(7.16)]{BGLV}, which control $\|\nabla v_{\mathrm{Ag}}(\cdot,t)\|_{L^\infty}$ and the $C^{1,\alpha}$ character of $\partial D_t$ on bounded time intervals. Indeed, $|\partial_s v_{\mathrm{Ag}}|\le\|\nabla v_{\mathrm{Ag}}\|_{L^\infty}$; by \eqref{CDE}, $|\partial_t g|\le\|\nabla v_{\mathrm{Ag}}\|_{L^\infty}\,|g|$ and $\bigl|\partial_t|\gamma(\xi,t)-\gamma(\eta,t)|\bigr|\le\|\nabla v_{\mathrm{Ag}}\|_{L^\infty}|\gamma(\xi,t)-\gamma(\eta,t)|$, which bound $|g|_*$ and $\Gamma(\gamma)$ by Gronwall's lemma; and $\|\tau(t)\|_{C^\alpha}$ is controlled by the $C^{1,\alpha}$ bound for $\partial D_t$ together with the bound for $g$. Thus, we can follow the computations in \cite[Proposition 3.3]{KL} to obtain the improved bound

\begin{equation}\pnorm{W^{2,p}(\mathbb{T})}{\vag(\gamma(t))}\le C_T(1+\pnorm{W^{2,p}(\T)}{\gamma(t)}), \qquad 0\leq t \leq T.
\label{improvedestimate}
\end{equation}

Since solutions $\gamma$ of equation \eqref{CDE} satisfy $\gamma(t)=\gamma_0+\int_0^t \vag(\gamma(s))\,\diff s$, using the bound in \eqref{improvedestimate} and Gronwall's lemma we get a $W^{2,p}$ bound for $\gamma(t)$ that does not blow up in finite time,
$$1+\pnorm{W^{2,p}(\mathbb{T})}{\gamma(t)}\le (1+\pnorm{W^{2,p}(\mathbb{T})}{\gamma_0})\exp(C_T t).$$

Thus, the maximum time of existence cannot be finite and the solution is global. Moreover, $v_{\mathrm{Ag}}(\gamma)\in C([0,\infty);W^{2,p}(\mathbb T))$ by the Lipschitz bounds above, so $\gamma\in C^1([0,\infty);W^{2,p}(\mathbb T))$, as in \cite[Theorem~3.4]{KL}.

\begin{theorem}[Global-in-time existence and uniqueness in $W^{2,p}$]
\label{thm:globalw2p}
    Let $1<p<\infty$. If $\gamma_0 \in X_p$ then there is a unique global-in-time solution $\gamma$ to equation \eqref{CDE} in $C([0,\infty);X_p)$  with initial condition $\gamma(0)=\gamma_0.$ Moreover, $\gamma\in C^1([0,\infty);W^{2,p}(\mathbb T))$.
\end{theorem}

\section{The curvature equation for aggregation}\label{sec:curvature}

In \cite{KL}, a differential equation for the curvature $\kappa$ is derived from geometric identities. This is equation (4.7) in their article and hence it is independent of the particular evolution equation
\begin{equation}
    \partial_t \kappa = -2\kappa \partial_s \partial_t \gamma \cdot \tau - \partial_s^2 \partial_t \gamma \cdot n.
    \label{evolcurv}
    \end{equation}

We start with the first term on the right-hand side of the equation. We need to give an expression for $\partial_s\partial_t \gamma \cdot \tau$. 

\begin{lemma}[First term in curvature equation]\label{integralmenosunmedio}
Let $\gamma\in C^{1,\alpha}(\T;\R^2)$, $0<\alpha<1$, be a regular
counterclockwise Jordan parametrization of the boundary of a bounded
domain $D$. Then, for every $\xi\in\T$,
\begin{equation}\label{eq:curv-half-flux}
 \frac{1}{2\pi}\int_{\T}
 \frac{n(\eta)\cdot(\gamma(\xi)-\gamma(\eta))}
 {|\gamma(\xi)-\gamma(\eta)|^2}\,g(\eta)\,\diff\eta
 =-\frac12.
\end{equation}
Consequently, for the aggregation velocity,
\begin{equation}\label{eq:curv-normal-integral}
 \frac{1}{2\pi}\int_{\T}
 [n(\eta)\cdot n(\xi)]
 \frac{(\gamma(\xi)-\gamma(\eta))\cdot n(\xi)}
 {|\gamma(\xi)-\gamma(\eta)|^2}\,g(\eta)\,\diff\eta
 =-\frac12-\partial_s\partial_t\gamma(\xi)\cdot\tau(\xi).
\end{equation}
\end{lemma}

\begin{proof}
Fix $\xi \in \T$ and set $x=\gamma(\xi)$ and $y=\gamma(\eta)$. Also, let $B_{\varepsilon}(x)$ denote the ball centered at $x$ with radius $\varepsilon$. It is easy to check that $\frac{x-y}{\left|x-y\right|^2}$ is a divergence free vector field. Then, by the divergence theorem applied to the domain $D\backslash B_{\varepsilon}(x)$ we have
$$\int_{\partial D\backslash B_{\varepsilon}(x)} \frac{(x-y)\cdot n(y)}{\left|x-y\right|^2}\,\diff\sigma(y)+\frac{\text{length}(D\cap \partial B_{\varepsilon}(x))}{\varepsilon}=0,$$
where we have used that the normal in the circular part is $(x-y)/\varepsilon.$ Taking the limit when $\varepsilon \to 0$ and putting the labels $x$ and $y$ as $\gamma(\xi)$ and $\gamma(\eta)$ we obtain the identity \eqref{eq:curv-half-flux}.

Projecting onto the normal and tangential directions, we have
\begin{align*}
    &\frac{1}{2\pi}\int_{\mathbb{T}}\frac{n(\eta)\cdot(\gamma(\xi)-\gamma(\eta))}{\left|\gamma(\xi)-\gamma(\eta)\right|^2}g(\eta)\,\diff\eta=\\&
    \hspace{0.5cm}=\frac{1}{2\pi}\int_{\mathbb{T}}\frac{[n(\eta)\cdot \tau(\xi)][(\gamma(\xi)-\gamma(\eta))\cdot\tau(\xi)]}{\left|\gamma(\xi)-\gamma(\eta)\right|^2}g(\eta)\,\diff\eta\\
    &\hspace{0.5cm}+\frac{1}{2\pi}\int_{\mathbb{T}}\frac{[n(\eta)\cdot n(\xi)][(\gamma(\xi)-\gamma(\eta))\cdot n(\xi)]}{\left|\gamma(\xi)-\gamma(\eta)\right|^2}g(\eta)\,\diff\eta.
\end{align*}
The left-hand side equals $-1/2$ by \eqref{eq:curv-half-flux} and the first term in the right-hand side equals $\partial_s \partial_t \gamma(\xi)\cdot \tau(\xi)$ by \eqref{firstderofv} and thus we have proved \eqref{eq:curv-normal-integral}.
\end{proof}
\begin{remark}
Although the geometric computations remain valid, rotating the velocity changes the projections entering the curvature equation. Then, the equivalent version of the previous lemma for Euler equation is not valid here.
\end{remark}

We take care now of the second term in the curvature equation.

\begin{lemma}[Second term in curvature equation]
Under the hypothesis in Lemma \ref{integralmenosunmedio}, if $\gamma \in W^{2,p}(\T;\mathbb{R}^2)$, $p>2$, we have 
\begin{equation}
    \partial_s^2\partial_t\gamma\cdot n = \left(\frac{1}{2}+\partial_s\partial_t \gamma\cdot \tau\right)\kappa -F, 
    \label{secondpartcurv}
    \end{equation}
where the remainder $F$ is locally integrable and can be written as
\begin{equation}\label{goodterms}
\begin{aligned}
 F(\xi):={}&-\frac{1}{2\pi}\int_{\T}
 \kappa(\eta)[\tau(\eta)\cdot n(\xi)]
 \frac{(\gamma(\xi)-\gamma(\eta))\cdot\tau(\xi)}
 {|\gamma(\xi)-\gamma(\eta)|^2}\,g(\eta)\,\diff\eta
 \\
 &-\frac{1}{2\pi}\int_{\T}[n(\eta)\cdot n(\xi)]
 \frac{(\tau(\xi)-\tau(\eta))\cdot\tau(\xi)}
 {|\gamma(\xi)-\gamma(\eta)|^2}\,g(\eta)\,\diff\eta
 \\
 &+\frac{1}{\pi}\int_{\T}[n(\eta)\cdot n(\xi)]
 \frac{\bigl[(\gamma(\xi)-\gamma(\eta))\cdot\tau(\xi)\bigr]
       \bigl[(\gamma(\xi)-\gamma(\eta))\cdot(\tau(\xi)-\tau(\eta))\bigr]}
 {|\gamma(\xi)-\gamma(\eta)|^4}\,g(\eta)\,\diff\eta.
\end{aligned}\end{equation}
\end{lemma}

\begin{proof}
Taking the normal projection of \eqref{secondderofv}, note that there is a term containing $\kappa(\xi)$, which is
$$-\frac{\kappa(\xi)}{2\pi}\int_{\T}
 [n(\eta)\cdot n(\xi)]
 \frac{(\gamma(\xi)-\gamma(\eta))\cdot n(\xi)}
 {|\gamma(\xi)-\gamma(\eta)|^2}\,g(\eta)\,\diff\eta.$$

 By \eqref{eq:curv-normal-integral} this is equal to
 $$\left(\frac{1}{2}+\partial_s\partial_t \gamma(\xi)\cdot \tau(\xi)\right)\kappa(\xi).$$
The rest of the terms of the normal projection of \eqref{secondderofv} are exactly $-F(\xi),$ which finishes the proof. 
\end{proof}
\begin{remark}
Note that the integral with principal value in \eqref{secondderofv} is integrable after taking the scalar product with $n(\xi)$ because
$$\left|\tau(\eta)\cdot n(\xi)\right| \le C \left|\eta-\xi\right|^{\alpha}.$$
This inequality is true for $\gamma \in C^{1,\alpha}(\mathbb{T}).$ See, for instance, \cite{KL}, Lemma 3.2.   
\end{remark}

\begin{remark}[Continuity of $F$]\label{rmk:cont}
Take $p=4$. By the preceding remark, the first integrand in \eqref{goodterms}
is bounded by $C|\kappa(\eta)|\,|\xi-\eta|^{-1/4}$.
The other two are bounded by $C|\xi-\eta|^{-1/2}$, using
$$  1-\tau(\xi)\cdot\tau(\eta)
  =\tfrac12|\tau(\xi)-\tau(\eta)|^2 $$
and the first-order Taylor estimate for $\gamma\in C^{1,3/4}$.
Since $\kappa\in C([0,T];L^4)$, Hölder's inequality bounds the
contribution from $|\xi-\eta|<\varepsilon$ by $C\varepsilon^{1/2}$,
uniformly in $\xi$ and $t$. The truncated integrals are continuous
in $(\xi,t)$, so their uniform limit satisfies
$$F\in C([0,T];C(\mathbb T)).$$
\end{remark}

Finally, by substituting identities \eqref{eq:curv-normal-integral} and \eqref{secondpartcurv} into \eqref{evolcurv} we obtain the desired equation for the curvature 

\begin{equation}\label{eq:curv-main}
 \partial_t\kappa
 =\left(-\frac12-3\,G\right)
  \kappa+F,
\end{equation}
where $F$ is defined in \eqref{goodterms} and $G(\xi):=\partial_s \partial_t \gamma(\xi)\cdot \tau(\xi)$.

Finding an equation describing the evolution of the metric $g$ is simpler. Since $g=\left|\partial_\xi \gamma\right|$, differentiating $g^2$ gives $\partial_t g = g \partial_s\vag \cdot \tau$. Again, the term $\partial_s \vag$ in \eqref{firstderofv} contains a singular integral but becomes absolutely integrable when multiplied by $\tau$. In fact,
\begin{equation}
\label{evolmetric}
\partial_t g(\xi)
=
\frac{g(\xi)}{2\pi}
\int_{\mathbb T}
\bigl[n(\eta)\cdot\tau(\xi)\bigr]\,
\frac{\bigl(\gamma(\xi)-\gamma(\eta)\bigr)\cdot\tau(\xi)}
     {|\gamma(\xi)-\gamma(\eta)|^2}\,
g(\eta)\,\diff\eta=g(\xi)G(\xi).
\end{equation}

Equations \eqref{eq:curv-main} and \eqref{evolmetric} establish a system of differential equations
\begin{equation}
    \label{systemkg}
    \tag{CM}
    \begin{cases}
        \partial_t \kappa = \left(-\frac12-3\,G\right)
  \kappa+F,\\
        \partial_t g=Gg,
    \end{cases}
\end{equation}
for $F$ defined in \eqref{goodterms} and $G$ defined in \eqref{evolmetric}. This is the curvature-metric system. Solving it is equivalent to solving \eqref{CDE}. The following statement is the analogue of {\cite[Proposition 4.2]{KL}} and it states this equivalence. We repeat it here for the reader's convenience, but omit the proof.

\begin{prop}[Equivalence CDE to curvature-metric equations]
Let $\gamma_0 \in W^{2,p}(\mathbb{T})$ be a proper parametrization
of a simple closed curve, and let
$(g_0,\kappa_0) \in W^{1,p}(\mathbb{T}) \times L^p(\mathbb{T})$
be the corresponding metric  and curvature.
Then the following statements hold.
\begin{enumerate}
    \item
    If $\gamma \in C([0,T];X_p)$ is a solution of \eqref{CDE}
    with initial data $\gamma_0$, then its metric $g$ and
    curvature $\kappa$ satisfy equations \eqref{systemkg}
    with initial data $g_0$ and $\kappa_0$.

    \item
    If
    $$        (g,\kappa) \in
        C\bigl([0,T];W^{1,p}(\mathbb{T}) \times L^p(\mathbb{T})\bigr),    $$
    regarded as the metric  and curvature, defines a simple
    closed curve and satisfies equations \eqref{systemkg}
    with initial data $g_0$ and $\kappa_0$, then there exists
    a solution
    $\gamma \in C([0,T];W^{2,p}(\mathbb{T}))$
    of \eqref{CDE} with initial data $\gamma_0$
    such that $(g,\kappa)$ is the metric and curvature
    of $\gamma$.
\end{enumerate}
\end{prop}

\section{\texorpdfstring{Persistence of $C^2$ regularity}{Persistence of C2}}\label{sec:C2agg}

In the previous section, we derived a system of differential equations involving the curvature and the metric \eqref{systemkg}. Now, we solve this system and verify that the solution provides a $C^2$ regular Jordan curve.

\persistence*

\begin{proof}

We can integrate the system \eqref{systemkg}. The term $G$ is well-defined and uniformly bounded in $L^{\infty}$ for every time interval. We can integrate the metric $g$ in terms of $G$ as
$$g(\xi,t)=g_0(\xi)\exp\left(\int_0^t G(\xi,s)\,\diff s\right).$$
Since $\left| G(\xi,t)\right|\le C_T$ for every $t\in[0,T]$ we have

$$g_0(\xi)\mathrm{e}^{-C_T t} \le g(\xi,t) \le g_0(\xi) \mathrm{e}^{C_T t}.$$

Thus, the metric remains uniformly bounded above and away from zero and so the parametrization remains regular for every time.

In order to achieve a closed formula for the curvature we differentiate with respect to time the function $t\to \mathrm{e}^{t/2}g(\xi,t)^3 \kappa(\xi,t)$, so $\mathrm{e}^{t/2} g^3$ is an integrating factor. Making use of the equation for $g$, we have
\begin{align*}
\partial_t(\mathrm{e}^{t/2}g^3 \kappa)=&\frac{1}{2}\mathrm{e}^{t/2}g^3\kappa + \mathrm{e}^{t/2}\left[3g^2\partial_t g\,\kappa + g^3 \partial_t \kappa\right]  \\
=&\frac{1}{2}\mathrm{e}^{t/2}g^3\kappa + \mathrm{e}^{t/2}\left[3g^3 G \kappa + g^3\left(\left(-\frac{1}{2}-3G\right)\kappa+F\right)\right] \\
=&\mathrm{e}^{t/2}g^3 F.
\end{align*}
Now integrating with respect to time between $0$ and $t$, we have
$$\kappa(\xi,t)=\frac{e^{-t/2}}{g(\xi,t)^3}\left[g_0(\xi)^3\kappa_0(\xi)+\int_0^t e^{u/2}g(\xi,u)^3F(\xi,u)\,du
\right].$$

By Remark \ref{rmk:cont} and the $ W^{2,4} $ regularity of
the boundary velocity, $ F $ and $ G $ are continuous on
$ \mathbb T\times[0,T] $. Since $ g>0 $ and $ \kappa_0 $ is continuous,
the formula above gives continuous curvature, so the arc-length
parametrization is $ C^2 $.

\end{proof}

\begin{remark}[$C^{1,1}$ patches and transport of curvature singularities]\label{rmk:C11}
The proof of Theorem~\ref{thm:C2persistence} only uses that $\gamma_0\in W^{2,4}$: Remark~\ref{rmk:cont} bounds $F$ in terms of $\|\kappa\|_{L^4}$, and Theorem~\ref{thm:globalw2p} provides these bounds for every $W^{2,4}$ datum. Hence, for every such datum and every $t\ge0$, the formula for $\kappa$ shows that
\[
\kappa(\cdot,t)-e^{-t/2}\Bigl(\frac{g_0}{g(\cdot,t)}\Bigr)^{3}\kappa_0
\]
has a continuous representative. Thus singularities of the initial curvature are only transported along the material labels and multiplied by a continuous positive factor; they are neither created nor amplified. In particular, if $\partial D_0$ is of class $C^{1,1}$, i.e.\ $\kappa_0\in L^\infty$, then $\partial D_t$ is of class $C^{1,1}$ for every $t\ge0$.
\end{remark}

\section{General kernel}\label{sec:generalK} 

Let $L:\R^2 \to \R^2$ be a linear transform. We consider now the transport equation with patch initial data
\begin{equation}
\label{generaltransport}
\begin{cases}
    \partial_t\rho+v_L\cdot \nabla \rho = 0,\\
    v_L(\cdot,t)=K_L\ast\rho(\cdot,t),\quad K_L=L(\nabla N),\\            
    \rho(\cdot,0)=\rho_0 = \chi_{D_0}.
    \end{cases}
\end{equation}

Thus, \eqref{generaltransport} reduces to 2D-Euler for $L=\begin{pmatrix}  0 & -1 \\ 1 & 0 \end{pmatrix}$ and becomes the aggregation transport equation for $L=-\operatorname{Id}.$ Hence, equation \eqref{generaltransport} is a generalization of both of them. 

A first remark is that, since $\kag=-\operatorname{Id}(\nabla N)$, then $v_L=-L\vag$, and hence some of the computations done until now can be easily adapted. On the other hand, the divergence of the velocity field can be more involved than in the previous scenarios. In the sense of distributions, we have
$$
\operatorname{div}v_L=\frac{\operatorname{tr}L}{2}\,\rho+T_L\rho,
\qquad
T_L\rho:=\operatorname{p.v.}\frac{(\operatorname{tr}L)\,|x|^2-2\,x\cdot Lx}{2\pi|x|^4}*\rho .
$$
For Euler ($L$ antisymmetric) both terms vanish, and for aggregation ($L=-\operatorname{Id}$) the kernel of $T_L$ vanishes, so that $\operatorname{div}\vag=-\rho=-\chi_D$. In general, $T_L$ is a Calderón--Zygmund operator with an even kernel, which vanishes only when the symmetric part of $L$ is a multiple of the identity. Hence, for a general $\rho\in L^\infty$, $\operatorname{div}v_L$ is only in $\operatorname{BMO}$. In the arguments below this difficulty in controlling the divergence is not relevant, but it is when dealing with the $C^{1,\alpha}$ persistence of regularity, which is an intermediate step of the proof. We use the main theorem in \cite{CMOV} to ensure this regularity.

Since $v_L=-L\vag$, integration by parts gives an associated Contour Dynamics Equation for the problem
\begin{equation}
    \label{LCDE}
    \tag{LCDE}
\partial_t\gamma(\xi,t)=-\frac{1}{2\pi}
\int_{\mathbb T}L n(\eta,t)\,\log|\gamma(\xi,t)-\gamma(\eta,t)|\,g(\eta,t)\,\diff\eta.
\end{equation}

The first thing we need is global-in-time existence and uniqueness of solutions for \eqref{LCDE}.

\begin{theorem}[Global-in-time existence and uniqueness in $W^{2,p}$ for \eqref{LCDE}]
\label{thm:lglobalw2p}
    Let $1<p<\infty$. If $\gamma_0 \in X_p$ then there is a unique global-in-time solution $\gamma$ to equation \eqref{LCDE} in $C([0,\infty);X_p)$  with initial condition $\gamma(0)=\gamma_0.$ Moreover, $\gamma\in C^1([0,\infty);W^{2,p}(\mathbb T))$.
\end{theorem}
\begin{proof}
The velocity in \eqref{LCDE} is the one considered in Section~\ref{sec:w2p} for a general constant matrix $L$, so the proofs of Theorems~\ref{thm:localw2p} and~\ref{thm:globalw2p} apply verbatim, with \eqref{eq:C1alpha} provided by \cite[Section~2]{CMOV}.
\end{proof}

Now we establish the curvature-metric system that will allow us to study the persistence of $C^2$ smoothness of the curve.

\begin{prop}[General curvature-metric system]\label{prop:Lsystem}
Let $p>2$. Along the solution $\gamma$ in Theorem \ref{thm:lglobalw2p} we have
\begin{equation}
    \label{Lsystem}
    \begin{cases}        \partial_tg=g\,\partial_s\partial_t\gamma\cdot\tau,\\
    \partial_t\kappa=c_L\kappa+q_LH\kappa+F_L,
  
    \end{cases}
\end{equation}
where 
$$q_L:=\frac12(L\tau)\cdot n,\quad\text{and}\quad c_L:={}-2\,\partial_s\partial_t\gamma\cdot\tau
 +(Ln\cdot n)\left(\frac12+\partial_s \vag\cdot\tau\right)-(L\tau\cdot n)(\partial_s \vag\cdot n)$$
 and $F_L$ is a sum of absolutely convergent integrals. Explicitly,
\begin{equation}\label{exprforF}
\begin{aligned}
 F_L(\xi)={}&\frac1{2\pi}\int_{\T}\kappa(\eta)
 \left\{
 [L\tau(\eta)\cdot n(\xi)]
 \frac{(\gamma(\xi)-\gamma(\eta))\cdot\tau(\xi)}
 {|\gamma(\xi)-\gamma(\eta)|^2}\,g(\eta)\right.\\
 &\hspace{36mm}\left.
 -\frac12[L\tau(\xi)\cdot n(\xi)]
       \cot\!\left(\frac{\xi-\eta}{2}\right)
 \right\}\diff \eta\\
 +&\frac1{2\pi}\int_{\T}[Ln(\eta)\cdot n(\xi)]
 \frac{(\tau(\xi)-\tau(\eta))\cdot\tau(\xi)}
 {|\gamma(\xi)-\gamma(\eta)|^2}\,g(\eta)\,\diff \eta\\
 -&\frac1\pi\int_{\T}[Ln(\eta)\cdot n(\xi)]
 \frac{[(\gamma(\xi)-\gamma(\eta))\cdot\tau(\xi)]
       [(\gamma(\xi)-\gamma(\eta))\cdot(\tau(\xi)-\tau(\eta))]}
 {|\gamma(\xi)-\gamma(\eta)|^4}\,g(\eta)\,\diff \eta.
\end{aligned}
\end{equation}
\end{prop}

 \begin{proof}
Throughout the proof, $\vag$ denotes the aggregation velocity
associated with the current domain, evaluated along its boundary.
Since $L$ is constant, it commutes  with time and spatial derivatives. Hence the identity for the evolution of the curvature yields
\begin{equation*}
\begin{aligned}
 \partial_t\kappa=-2\kappa\,\partial_s\partial_t\gamma\cdot\tau
   -\partial_s^2\partial_t\gamma\cdot n=-2\kappa\,\partial_s\partial_t\gamma\cdot\tau
   +L(\partial_s^2\vag)\cdot n.
\end{aligned}
\end{equation*}
We therefore use the formula for $\partial_s^2\vag$
obtained in the aggregation case, apply $L$, and project onto $n$.
The two terms in that formula which do not explicitly contain
curvature become precisely the second and third integrals in
\eqref{exprforF}. It remains to identify the contributions of the
terms containing $\kappa(\xi)$ and $\kappa(\eta)$.

The term containing $\kappa(\xi)$ is
$$ -\frac{\kappa(\xi)}{2\pi}
 \int_{\T}[Ln(\eta)\cdot n(\xi)]
 \frac{(\gamma(\xi)-\gamma(\eta))\cdot n(\xi)}
      {|\gamma(\xi)-\gamma(\eta)|^2}
 \,g(\eta)\,\diff \eta.$$
To evaluate it, we need both components, in the frame at $\xi$,
of the vector integral before applying $L$.
Its normal component is the identity \eqref{eq:curv-normal-integral}:
\begin{equation*}
\begin{aligned}
 &\frac1{2\pi}\int_{\T}[n(\eta)\cdot n(\xi)]
 \frac{(\gamma(\xi)-\gamma(\eta))\cdot n(\xi)}
      {|\gamma(\xi)-\gamma(\eta)|^2}
 \,g(\eta)\,\diff \eta\\
 &\hspace{25mm}
 =-\frac12-\partial_s \vag(\xi)\cdot\tau(\xi).
\end{aligned}
\end{equation*}
The tangential component simply satisfies
\begin{equation*}
\begin{aligned}
 &\frac1{2\pi}\int_{\T}[n(\eta)\cdot\tau(\xi)]
 \frac{(\gamma(\xi)-\gamma(\eta))\cdot n(\xi)}
      {|\gamma(\xi)-\gamma(\eta)|^2}
 \,g(\eta)\,\diff \eta\\
 &\hspace{25mm}
 =\partial_s \vag(\xi)\cdot n(\xi).
\end{aligned}
\end{equation*}

Indeed, subtracting from the left-hand side the normal component of \eqref{firstderofv}, the numerator becomes
$$[n(\eta)\cdot\tau(\xi)][(\gamma(\xi)-\gamma(\eta))\cdot n(\xi)]-[n(\eta)\cdot n(\xi)][(\gamma(\xi)-\gamma(\eta))\cdot\tau(\xi)]
=\tau(\eta)\cdot(\gamma(\eta)-\gamma(\xi)),$$
as one sees by writing $n(\eta)$ and $\gamma(\xi)-\gamma(\eta)$ in the orthonormal frame $(\tau(\xi),n(\xi))$ and using $\tau^\perp=-n$, $n^\perp=\tau$. After multiplication by $g(\eta)/|\gamma(\xi)-\gamma(\eta)|^2$ this is $\partial_\eta\log|\gamma(\xi)-\gamma(\eta)|$, whose principal value integral over $\T$ vanishes.

Now, decomposing $n(\eta)$ in the same frame and using the identities for the normal and tangential components, the term containing $\kappa(\xi)$ becomes
\[
 \kappa(\xi)\left\{
 [Ln(\xi)\cdot n(\xi)]
 \left(\frac12+\partial_s \vag(\xi)\cdot\tau(\xi)\right)
 -
 [L\tau(\xi)\cdot n(\xi)]
 [\partial_s \vag(\xi)\cdot n(\xi)]
 \right\}.
\]
Together with the geometric term
$-2\kappa\,\partial_s\partial_t\gamma\cdot\tau$,
this gives exactly $c_L\kappa$.

We then focus on the term containing $\kappa(\eta)$
\begin{equation}\frac1{2\pi}\operatorname{p.v.}\int_{\T}
 \kappa(\eta)[L\tau(\eta)\cdot n(\xi)]
 \frac{(\gamma(\xi)-\gamma(\eta))\cdot\tau(\xi)}
      {|\gamma(\xi)-\gamma(\eta)|^2}
 \,g(\eta)\,\diff \eta.
 \label{Lhilterterm}
\end{equation}
This will be compared with a Hilbert transform as in Kiselev-Luo. In fact, we will follow exactly the arguments in  \cite[Section 5.1]{KL}.
Specifically, in the singular integral containing $\kappa(\eta)$ we add and subtract
$q_L(\xi)\mathcal{H}\kappa(\xi)$, where
$q_L(\xi)=\frac12L\tau(\xi)\cdot n(\xi)$ and $\mathcal{H}$ is the Hilbert transform.
At the level of the kernel, this corresponds to adding and subtracting
\begin{equation}
 \frac12[L\tau(\xi)\cdot n(\xi)]
 \cot\!\left(\frac{\xi-\eta}{2}\right).
 \label{addsubstractL}
\end{equation}
By the same diagonal comparison and periodization used in
\cite[Proposition 5.2]{KL}, this kernel
matches the leading singularity of the original one.
The added term contributes exactly $q_L\mathcal{H}\kappa$, while subtracting  \eqref{addsubstractL} from the kernel in \eqref{Lhilterterm}  and integrating again against $\kappa(\eta)\diff \eta$ gives the first integral in \eqref{exprforF} and is
absorbed into $F_L$.
The cancellation makes this difference absolutely integrable
for $p>2$. The other two integrals in \eqref{exprforF} are the
remaining terms of the curvature equation.

For $p=4$, the bounds of Remark \ref{rmk:cont} also hold
for the integrands in \eqref{exprforF}. Indeed, near the diagonal, multiplying
the kernel in braces by $ \xi-\eta $ gives a $ C^{3/4} $ function
of $ \xi $ vanishing at $ \xi=\eta $, since
$ (\gamma(\xi)-\gamma(\eta))/(\xi-\eta) $ is an average of $ \gamma' $.
For $ |\xi-\eta|\ge 2|\delta| $, taking differences in the chord and
tangent factors gives the same bounds for the increment from $ \xi $
to $ \xi+\delta $, multiplied by $ (|\delta|/|\xi-\eta|)^{3/4} $.
Splitting at $ |\xi-\eta|=2|\delta| $ and using $ \kappa\in L^4 $ yields
$$
\sup_{0\le t\le T}\|F_L(\cdot,t)\|_{C^{1/2}(\mathbb T)}<\infty.
$$

\end{proof}

\begin{remark}
If $L=c\operatorname{Id}$, the Hilbert transform term in the previous proposition vanishes; in particular, for $L=-\operatorname{Id}$ we recover the curvature-metric system \eqref{systemkg}. The converse is transparent in complex notation. Identifying $\R^2$ with $\mathbb C$, every linear map of $\R^2$ can be written as $Lw=\alpha w+\beta\bar w$ with $\alpha,\beta\in\mathbb C$, and then $K_L(z)=\frac{\alpha}{2\pi\bar z}+\frac{\beta}{2\pi z}$. Writing $\tau=e^{i\theta}$, so that $n=-i\tau$, we get
$$q_L=-\tfrac12\bigl(\operatorname{Im}\alpha+\operatorname{Im}(\beta e^{-2i\theta})\bigr).$$
Moreover, $q_L$ vanishes for every direction $\theta$ if and only if $\alpha\in\R$ and $\beta=0$, that is, $L=c\operatorname{Id}$. Hence, if $L$ is not a multiple of the identity, there is a direction $\tau$ with $q_L\neq0$.
\end{remark}

Now we are prepared to characterize the linear maps $L$ for which $C^2$ boundary regularity persists.

\classification*

\begin{proof}
If we apply the Hilbert transform $\mathcal{H}$ to the second equation in \eqref{Lsystem}, we get
\begin{align*}
    \partial_t \mathcal{H} \kappa &= \mathcal{H}c_L\kappa + \mathcal{H}q_L\mathcal{H}\kappa + \mathcal{H}F_L \\
    &=c_L \mathcal{H}\kappa + [\mathcal{H},c_L]\kappa + q_L \mathcal{H}^2 \kappa + [\mathcal{H},q_L]\mathcal{H}\kappa + \mathcal{H}F_L \\
    &=c_L \mathcal{H}\kappa -q_L \kappa + \tilde{F}_L,
\end{align*}
where $[T_1,T_2]f=T_1(T_2f)-T_2(T_1f)$ is a commutator and where we have defined 
$$\tilde{F}_L:= \mathcal{H}F_L+[\mathcal{H},c_L]\kappa + [\mathcal{H},q_L]\mathcal{H}\kappa+q_L \langle \kappa \rangle,$$

where $\langle f \rangle := \fint_{\mathbb{T}} f(\eta)\,\diff \eta.$
Also, we have used $\mathcal{H}^2 f=-f+ \langle f \rangle $. Hence we can write a system of differential equations for $\kappa$ and $\mathcal{H}\kappa$ as
\begin{equation}
    \label{systemkappahilber}
\partial_t\begin{pmatrix} \kappa \\ \mathcal{H}\kappa \end{pmatrix} = \begin{pmatrix} c_L & q_L \\ -q_L & c_L \end{pmatrix} \begin{pmatrix} \kappa \\ \mathcal{H}\kappa \end{pmatrix}+\begin{pmatrix} F_L \\ \tilde{F}_L \end{pmatrix}
\end{equation}

Since $\mathcal{H}$ is bounded on $C^{1/2}$, the estimate in Proposition \ref{prop:Lsystem}
also bounds $ \mathcal{H}F_L $ in $ C^{1/2} $. Lemma~5.5 of~\cite{KL}, with $ p=4 $
and $ \sigma=3/4 $, gives the same bound for both commutators.

We can solve the system of differential equations \eqref{systemkappahilber} by solving the homogeneous part and then by variation of constants. For each time $t$, the following identity holds for almost every $\xi$.
\begin{equation}
\begin{pmatrix}
    \kappa(\xi,t)\\
    \mathcal{H}\kappa(\xi,t)
\end{pmatrix}
=
a(\xi,t)
\begin{pmatrix}
    \cos b(\xi,t)&\sin b(\xi,t)\\
    -\sin b(\xi,t)&\cos b(\xi,t)
\end{pmatrix}
\begin{pmatrix}
    \kappa_0(\xi)\\
    \mathcal{H}\kappa_0(\xi)
\end{pmatrix}
+
\begin{pmatrix}
    R(\xi,t)\\
    \tilde{R}(\xi,t)
\end{pmatrix},
\end{equation}
where
$$a(\xi,t)=\exp\!\left(\int_0^t c_L(\xi,s)\,\diff s\right),
\qquad
b(\xi,t)=\int_0^t q_L(\xi,s)\,\diff s.
$$
and
$$\begin{pmatrix}
    R(\xi,t)\\
    \tilde{R}(\xi,t)
\end{pmatrix}
=
\int_0^t \frac{a(\xi,t)}{a(\xi,s)}
\begin{pmatrix}
    \cos\bigl(b(\xi,t)-b(\xi,s)\bigr)&\sin\bigl(b(\xi,t)-b(\xi,s)\bigr)\\
    -\sin\bigl(b(\xi,t)-b(\xi,s)\bigr)&\cos\bigl(b(\xi,t)-b(\xi,s)\bigr)
\end{pmatrix}
\begin{pmatrix}
    F_L(\xi,s)\\
    \widetilde{F}_L(\xi,s)
\end{pmatrix}
\,\diff s.
$$
 Thus $R$ and $\widetilde R$ are uniformly bounded in $C^{1/2}$,
since the coefficients in the time integral are uniformly $C^{3/4}$
in $\xi$.  In particular, for a given time $t>0$
$$
\kappa(t)
=a(t)\cos b(t)\,\kappa_0
+a(t)\sin b(t)\,\mathcal{H}\kappa_0
+R(t).
$$

If $L=c\operatorname{Id}$ then $q_L=0$ so there is no Hilbert contribution in $\kappa(t)$ and we can proceed as in the aggregation case to show persistence of $C^2$ regularity of the boundary.

On the other hand, if $L\neq c\operatorname{Id}$, we choose an initial curve $\gamma_0$ as constructed in \cite[Equation (5.25)]{KL} orientated such that $q_L(0,0)\neq 0$. We can do such a construction since there exists a direction $\tau$ such that 
$$q_L=\frac{1}{2}(L\tau)\cdot n \neq 0$$
provided $L$ is not a multiple of the identity. There is a neighborhood $U$ of $(\xi=0,t=0)$ such that $\operatorname{sign}(q_L(\xi,t))=\operatorname{sign}(q_L(0,0)).$ Consequently, for $(\xi,t)\in U$ we have ensured $\left|\sin(b(\xi,t))\right|>0$ for $t>0$. 

Then, for the bounded curvature $\kappa_0$ associated to $\gamma_0$ we have that $\kappa(t)\notin L^{\infty}(\T)$ for short time, which proves the theorem.
\end{proof}
\begin{remark}[$C^{1,1}$ version of Theorem~\ref{thm:C2-classification}]\label{rmk:C11L}
Theorem~\ref{thm:C2-classification} also holds with $C^{1,1}$ in place of $C^2$. If $L=c\operatorname{Id}$, then $e^{-ct/2}g^3$ is an integrating factor for \eqref{Lsystem}, and as in Remark~\ref{rmk:C11},
\[
\kappa(\cdot,t)-e^{ct/2}\Bigl(\frac{g_0}{g(\cdot,t)}\Bigr)^{3}\kappa_0
\]
has a continuous representative, so $C^{1,1}$ regularity persists. If $L\neq c\operatorname{Id}$, the initial domain used in the proof is of class $C^2\subset C^{1,1}$ and its curvature becomes unbounded, so $C^{1,1}$ regularity is lost instantaneously.
\end{remark}

\section*{Funding and Acknowledgments}
    This research is supported by the project ``Análisis y Ecuaciones en Derivadas Parciales'' (Ref. PID2024-155320NB-I00) of the Spanish Ministry of Science, Innovation and Universities.

    We thank Joan Verdera for his helpful comments on the introduction.

\section*{Declaration on the use of AI}
The results of this paper and their original proofs were obtained by the authors without the use of artificial intelligence (AI) tools. Afterwards, we used GPT-6 Astra (OpenAI) and Claude Opus 5.5 (Anthropic) to improve the writing and to check the manuscript. This included verifying the steps of the proofs and suggesting corrections and more detailed justifications of some of them, which the authors checked before including them. The authors have reviewed and edited all the content and take full responsibility for this publication.

\bibliographystyle{amsplain}
\bibliography{biblio}

\end{document}